\documentstyle[12pt,leqno]{article}
\begin{document}
\title{Connections under Symplectic Reduction}
\author{{\normalsize by}\\Izu Vaisman}
\date{}
\maketitle
{\def\thefootnote{*}\footnotetext[1]%
{{\it 2000 Mathematics Subject Classification}
53D20. \newline\indent{\it Key words and phrases}:
Symplectic Reduction. Symplectic Connection. Presymplectic Connection.}}
\centerline{\footnotesize Dedicated to the memory of Prof. Gheorghe Vr\u
anceanu}
\begin{center} \begin{minipage}{12cm}
A{\footnotesize BSTRACT. In this note,
we give conditions which ensure the reduction of a
symplectic connection in the process of a Marsden-Weinstein reduction and of
the reduction of a presymplectic manifold.}
\end{minipage} \end{center}
\vspace{5mm}\noindent
Symplectic reduction is a technique which produces new symplectic manifolds
from either symplectic manifolds with symmetries or presymplectic manifolds
\cite{LM}. The aim of this note is to formulate conditions which ensure
that a compatible connection on the initial manifold induces a compatible
connection on the reduced manifold. Everything in this note is of the
class $C^{\infty}$.

For me, this subject is related with the name of Vr\u anceanu because my
Ph.D. thesis \cite{teza} was on the subject of symplectic connections, and
Vr\u anceanu was one of the referees.
\section{Marsden-Weinstein Reduction}
We refer the reader to \cite{LM} for the original Marsden-Weinstein
reduction theorem. In this theorem, the basic configuration consists of:
(i) a symplectic manifold $M^{2n}$ with the symplectic form $\omega$;
(ii) a Hamiltonian action $\Phi:G\times M\rightarrow M$ of a Lie group $G$
on $M$ with an equivariant moment map $J:M\rightarrow{\cal G}^{*}$ (the
dual space of the Lie algebra ${\cal G}$ of $G$);
(iii) the level set $C:=J^{-1}(\xi)$ of a non critical value $\xi\in{\cal
G}^{*}$.
Then, if $G_{\xi}$ is the isotropy subgroup of $\xi$ for the coadjoint
action of $G$ on ${\cal G}^{*}$, the group $G_{\xi}$ acts on $C$,
and, if $N:=C/G_{\xi}$ is a manifold,
$\omega$ induces a symplectic structure $\omega'$ on $N$.

Let us assume that $(M,\omega)$ is also equipped with a {\em symplectic
connection} $\nabla$ i.e., a torsionless linear connection on $M$ such that
$\nabla\omega=0$
(e.g., \cite{V2}). In order to formulate a corresponding reduction theorem,
we need the following notions: a) if the action $\Phi$ preserves the
connection $\nabla$ we say that $\Phi$ is an {\em affine action} \cite{KN},
b) a submanifold $C\subseteq M$ such that $TC$
is preserved by $\nabla$-parallel
translations along paths in $C$ is called a {\em self-parallel
submanifold}. Then, we get
\proclaim 1.1 Theorem. Let $(M,\omega)$ be a symplectic manifold with a
symplectic connection $\nabla$ and an affine Hamiltonian action $\Phi$ of
the Lie group $G$ which has the equivariant moment map $J$. Let
$\xi\in{\cal G}^{*}$ be a non critical value of $J$ such that $C:=J^{-1}(\xi)$
is $\nabla$-self-parallel, and such that the reduced symplectic manifold
$(N=C/G_{\xi},\omega')$ exists. Then $\nabla$ induces a well defined
symplectic connection $\nabla'$ on $(N,\omega')$. \par
\noindent{\bf Proof.} Since $C$ is self-parallel, the values of
$\nabla_{X}Y$ are tangent to $C$ whenever $X,Y\in \Gamma TC$ ($\Gamma$
denotes spaces of global cross sections of bundles), and $\nabla$ may also
be seen as a linear connection on $C$.

On the other hand, the connected components of the orbits of $G_{\xi}$
in $C$ define an isotropic foliation ${\cal I}$ such that $T{\cal I}=
(TC)^{\perp_{\omega}}\cap TC$, and, since $\nabla$ is a
symplectic connection, the bundle $T{\cal I}$ is parallel with respect to
$\nabla$.

Now, we will use the connection $\nabla$ to define an
induced connection $\nabla'$ on $N$ as follows. If ${\cal X},{\cal Y}
\in\Gamma TN$, we have ${\cal X}=\pi_{*}X$, ${\cal Y}=\pi_{*}Y$, where
$\pi:C\rightarrow C/G_{\xi}$ is the natural projection and $X,Y$ are
$G_{\xi}$-invariant vector fields on $C$ defined up to a term in $\Gamma
T{\cal I}$. Accordingly, we shall take
$$\nabla'_{{\cal X}}{\cal Y}=\pi_{*}(\nabla_{X}Y),
\leqno{(1.1)}$$ and show that
the result is well defined. The reasons for that are:\\
\indent1) the result of (1.1) does not change if
$Y\mapsto Y+Z$, $Z\in\Gamma T{\cal I}$, because $T{\cal I}$ is
$\nabla$-parallel;\\
\indent2) the same is true if $X\mapsto X+Z$, $Z\in\Gamma T{\cal I}$, since
$$\nabla_{Z}Y=\nabla_{Y}Z+[Y,Z]\in\Gamma T{\cal I};\leqno{(1.2)}$$
(1.2) is justified by the fact that
$\nabla$ has no torsion and preserves
$T{\cal I}$, and by the fact that,
locally, $Z=\sum_{i}\varphi_{i}Z_{i}$ where
$\varphi_{i}\in C^{\infty}(C)$ and $Z_{i}\in\Gamma T{\cal I}$
are infinitesimal actions of
elements of the Lie algebra ${\cal G}_{\xi}$ on $C$, therefore,
$[Y,Z]=\sum_{i}(Y\varphi_{i})Z_{i}\in\Gamma T{\cal I}$ because of the
invariance of $Y$;\\
\indent3) in (1.1), the vector field $\nabla_{X}Y$
is projectable since the action $\Phi$ of $G$ is affine, which
means that, $\forall g\in G$, the transformation
$\Phi_{g}(x):=\Phi(g,x)$ ($x\in M$) satisfies the condition
$$\Phi_{g*}(\nabla_{X}Y)=\nabla_{\Phi_{g*}X}
(\Phi_{g*}Y).\leqno{(1.3)}$$
Q.e.d.

It would be interesting to have conditions which imply the fact that the
level submanifold $C$ of Theorem 1.1 is self-parallel. We can give one such
condition which, unfortunately, is not simple to use. Namely, let us
agree to say that the moment map $J$ is {\em affine} if, whenever $J$ is
constant along a path $\gamma$ in $M$, the differential $J_{*}$ is constant
along any $\nabla$-horizontal lift $\tilde\gamma$ of $\gamma$ to $TM$. It is
easy to see that if $J$ is an affine moment map, $C=J^{-1}(\xi)$ is
self-parallel for all the non critical values
$\xi\in{\cal G}^{*}$. Indeed,
for $Y\in\Gamma TM$ with $Y/_{C}\in\Gamma TC$, and if $X_{0}\in T_{x_{0}}C$
$(x_{0}\in C)$, one has \cite{KN}
$$(\nabla_{X_{0}}Y)_{x_{0}}=\lim_{t\rightarrow 0}\frac{1}{t}
[\tau^{t}_{0}(Y(x(t)))-Y(x_{0})],$$
where $x(t)$ is a path in $M$ such that $x(0)=x_{0}$,
$\stackrel{.}{x}(0)=X_{0}$, and $\tau^{t}_{0}: T_{x(t)}M\rightarrow
T_{x_{0}}M$ is the $\nabla$-parallel translation along $x(t)$.
Accordingly, if $x(t)$ is in $C$ then, because of the affine character of
$J$, $(\nabla_{X_{0}}Y)_{x_{0}}\in T_{x_{0}}C$, and the
submanifold $C$ is self-parallel.

An example where Theorem 1.1 applies will be given at the end of
the next section.
\section{Cotangent Bundles}
An important example of a Hamiltonian action as required in reduction
theory is provided by the natural lift $\Phi^{*}$ of an action $\Phi$ of a
Lie group $G$ on a manifold $P$ to the cotangent bundle $T^{*}P$.
The latter is endowed with the canonical symplectic structure
$\omega=-d\lambda$, where $\lambda$ is the {\em Liouville $1$-form}
$$\lambda_{\alpha}(X)=<\alpha,\pi_{*}X>,\hspace{1cm}\alpha\in T^{*}P,\,
X\in T_{\alpha}T^{*}P,\,\pi:T^{*}P\rightarrow P,\leqno{(2.1)}$$
and $\Phi^{*}$ has the moment map $J:T^{*}P\rightarrow {\cal G}^{*}$ given
by
$$<J(\alpha),A>=<\alpha,A_{P}(\pi(\alpha))>,\leqno{(2.2)}$$
where $\alpha$ is as in (2.1), $A\in {\cal G}$,
the Lie algebra of $G$, and $A_{P}$ is the
infinitesimal transformation defined by $A$ on $P$ \cite{LM}.

In this section, we assume that the action $\Phi$ is affine with respect to
a torsionless connection $\nabla^{P}$ on $P$, and show how to get a symplectic
connection $\nabla$ on $T^{*}P$ with respect to which the action $\Phi^{*}$
is affine. This sets the scene for possible applications of Theorem 1.1.

Let ${\cal V}$ be the {\em vertical distribution} tangent to the fibers of
$M:=T^{*}P$, and
${\cal H}$ be the horizontal distribution of the connection
$\nabla^{P}$ on $M$. The following formula defines a connection
$\nabla^{\cal H}$ on the vector bundle ${\cal H}$:
$$\nabla^{\cal H}_{X}Y=\left\{ \begin{array}{lll}
pr_{\cal H}[X,Y]&{\rm if}&X\in\Gamma{\cal V},\,Y\in\Gamma{\cal
H},\vspace{1mm}\\
\widetilde{\nabla^{P}_{\pi_{*}X}\pi_{*}Y}&{\rm if}&
X\in\Gamma_{pr}{\cal H},\,Y\in\Gamma_{pr}{\cal H},
\end{array}\right. \leqno{(2.3)}$$
where tilde denotes the operation of taking the horizontal lift,
$pr$ stands for {\em projectable}, and projectability is with
respect to $\pi$. If $Y$ of the second line of (2.3) is of the general form
$\sum_{i}\varphi_{i}Y_{i}$, where $\varphi_{i}\in C^{\infty}(M)$,
$Y_{i}\in\Gamma_{pr}{\cal H}$,
$\nabla^{\cal H}_{X}Y$ is to be derived from (2.3) by means of the general
properties of a covariant derivative. Notice also that the first line of
(2.3) is equivalent to $\nabla^{\cal H}_{X}Y=0$ for all
$X\in\Gamma{\cal V}$ and for all $Y\in\Gamma_{pr}{\cal H}$.
It is easy to understand that (2.3) indeed provides a connection on ${\cal
H}$; in fact, it is a Bott connection with respect to the foliation by the
fibers of $T^{*}P$ \cite{Ml}.

The connection $\nabla^{\cal H}$ transposes to a connection
$\nabla^{{\cal H}^{*}}$ on the dual bundle ${\cal H}^{*}$ of ${\cal H}$.
On the other hand, the {\em musical morphism}
$\sharp_{\omega}:=\flat_{\omega}^{-1}$ is an
isomorphism between ${\cal H}^{*}$ and ${\cal V}$ (e.g., \cite{V3}). Hence,
$\nabla^{\cal V}:=\sharp_{\omega}\circ\nabla^{{\cal
H}^{*}}\circ\flat_{\omega}$ is a well defined connection on the vector
bundle ${\cal V}$.

Accordingly, we obtain a natural lift $\nabla^{M}:=\nabla^{\cal H}\oplus
\nabla^{\cal V}$ of the connection $\nabla^{P}$ to a connection on
$T^{*}P$. This connection satisfies the condition $\nabla^{M}\omega=0$.
Indeed, let $(x^{i})$ be local coordinates on $P$ such that
$$\nabla^{P}_{\frac{\partial}{\partial
x^i}}\frac{\partial}{\partial x^{j}}
=\Gamma^{k}_{ji}\frac{\partial}{\partial x^k}\hspace{1cm}
(\Gamma^{k}_{ji}=\Gamma^{k}_{ij}),\leqno{(2.4)}$$
and let $(y_{i})$ be the corresponding covector coordinates. Then,
$(x^{i},y_{i})$ are local coordinates on $M$, and \cite{V3}
$${\cal V}=span\left\{\frac{\partial}{\partial y_i}\right\},\;
{\cal H}= span\left\{X_{i}:=\frac{\partial}{\partial x^i}+
\Gamma^{s}_{ik}y_{s}\frac{\partial}{\partial y_k}\right\},
\leqno{(2.5)}$$
$${\cal V}^{*}=span\{\theta_i:=dy_i-\Gamma^{s}_{ik}y_sdx^k\},\;
{\cal H}^{*}=span\{dx^i\}, \leqno{(2.6)}$$
$$\omega=dx^i\wedge dy_i=dx^i\wedge\theta_i. \leqno{(2.7)}$$
In all these formulas the Einstein summation convention is used.

Accordingly, we get the following local equations of $\nabla^M$
$$\begin{array}{ll}
\nabla^{M}_{\frac{\partial}{\partial y_i}}X_j=0,&
\nabla^{M}_{X_{i}}X_{j}=\;\;\;\Gamma^{k}_{ji}X_{k},\vspace{1mm}\\
\nabla^{M}_{\frac{\partial}{\partial y_i}}\frac{\partial}{\partial y_j}=0,&
\nabla^{M}_{X_{i}}\frac{\partial}{\partial y_j}=
-\Gamma^{j}_{ik}\frac{\partial}{\partial y_{k}},
\end{array} \leqno{(2.8)}$$
and $\nabla^M\omega=0$ follows by straightforward computations.

Furthermore, if we come back to the $\nabla^P$-affine action $\Phi$, its
affine character implies the invariance of
the horizontal distribution ${\cal H}$ by $\Phi^*$. Then,
if we use $Y\in\Gamma_{pr}{\cal H}$ on the first line of (2.3) it is easy
to check that the connection $\nabla ^M$ satisfies condition (1.3).
Therefore, $M$ is endowed with the
$\omega$-compatible connection $\nabla^M$ and the affine
Hamiltonian action $\Phi^*$.

This is not yet the required situation since $\nabla^M$ may have torsion
\cite{V3}.
To fix this, let
$$^t\nabla^{M}_{X}Y:=\nabla^{M}_{Y}X+[X,Y]$$
be the {\em transposed connection}, and
$$\stackrel{\circ}{\nabla}:=\frac{1}{2}(\nabla_{X}Y+\,^t\nabla_{X}Y)$$
be the connection on $M$ known as the {\em symmetric part} of the
connection $\nabla^M$. Then $\stackrel{\circ}{\nabla}$ is a torsionless
connection which satisfies the condition (1.3), hence, it is preserved by
the action $\Phi^*$ but, it is no more symplectic.

However, there exists a well known transformation of such a connection
which yields a symplectic one. Namely (e.g., \cite{V2}):
$$\nabla_{X}Y=\stackrel{\circ}{\nabla}_{X}Y+A(X,Y),\leqno{(2.9)}$$
where $A(X,Y)$ is defined by the relation
$$\omega(A(X,Y),Z)=\frac{1}{2}(\stackrel{\circ}{\nabla}_X\omega)(Y,Z)
\leqno{(2.10)}$$
$$+\frac{1}{6}\{(\stackrel{\circ}{\nabla}_Y\omega)(X,Z)
+(\stackrel{\circ}{\nabla}_Z\omega)(X,Y)\}.$$
Since both $\stackrel{\circ}{\nabla}$ and $\omega$ are invariant by the
action $\Phi^*$ of $G$, the connection $\nabla$ of (2.9) satisfies
condition (1.3). Therefore, the action $\Phi^*$ is affine with respect to
the connection $\nabla$, and we are done.

Now, as a corollary of the above construction and of Theorem 1.1, we have
\proclaim 2.1 Proposition. Let $P$ be an arbitrary manifold with a linear
connection $\nabla^P$, and an affine action $\Phi$ of a Lie group $G$. Then,
the cotangent bundle $T^*P$ with the canonical symplectic structure
$\omega$ has a well defined symplectic connection $\nabla$
(the lift of $\nabla^{P}$) which is preserved
by the lift $\Phi^*$ of $\Phi$.
Furthermore, if there exists $\xi\in
{\cal G}^*$ which is non critical for the natural moment map $J$ of
$\Phi^*$ (given by (2.2)),
and such that $J^{-1}(\xi)$ is $\nabla$-self-parallel, the
connection $\nabla$ reduces to a symplectic connection of the reduced
manifold $J^{-1}(\xi)/G_{\xi}$ (if the latter exists).\par
We end this section by a simple example which illustrates both Proposition
2.1 and Theorem 1.1. Take $P={\bf R}^{n}$ with the natural flat torsionless
connection $\nabla$. Then, any affine subgroup on ${\bf R}^{n}$ is affine
with respect to $\nabla$, and we will use the group $G$ given by
$$\tilde x^{a}=sx^{a}+t^{a},\;\tilde x^{u}=x^{u},\leqno{(2.11)}$$
where $a=1,\ldots, h\leq n$, $u=h+1,\ldots n$,
$(x^{a},x^{u})$ are natural coordinates in ${\bf R}^{n}$, and $s,t^{a}\in{\bf
R},\,s\neq0$.

Then, $T^{*}P={\bf R}^{2n}$ with coordinates $(x^{i},y_{i})$
$(i=1,...,n)$, and with the
canonical symplectic form (2.7). Furthermore, in (2.4)
$\Gamma^{k}_{ji}=0$, and it follows from (2.8) that the lifted connection
$\nabla$ is exactly the
natural flat torsionless connection of ${\bf R}^{2n}$.
The action defined by (2.11) on covector coordinates
$(y_{i})$
is given by
$$\tilde y_{a}=\frac{1}{s}y_{a},\,\tilde y_{u}=y_{u}.\leqno{(2.12)}$$
Thus, the resulting action on ${\bf R}^{2n}$ is
affine with respect to the flat connection $\nabla$.

An infinitesimal action of the Lie algebra ${\cal G}$ of $G$ is obtained
from (2.11) by taking $s=1+\epsilon\sigma,\,t^{a}=\epsilon\tau^{a}$,
where $(\sigma,\tau^{a})$ is a vector in
${\cal G}$, then,
taking the derivative with respect to $\epsilon$ at $\epsilon=0$.

Now, (2.2) shows that the moment map $J$ has the expression
$$J(x,y)=(\sum_{a=1}^{h}x^{a}y_{a},y_{1},...,y_{h})\in{\cal
G}^*,\leqno{(2.13)}$$ and, for instance,
it follows that $\xi=(0,1,...,1)\in{\cal G}^{*}$ is a non critical
value of $J$. The corresponding level set $J^{-1}(\xi)$ has the equations
$$y_{1}=1,...,y_{h}=1,\sum_{a=1}^{h}x^{a}=0,\leqno{(2.14)}$$
and it is a plane of dimension $2n-h-1$ in ${\bf R}^{2n}$. Therefore,
$J^{-1}(\xi)$ is self-parallel with respect to the flat connection
$\nabla$. From (2.12), (2.13), we see that the isotropy subgroup $G_{\xi}$
consists of the transformations (2.11) where $s=1$ i.e., translations of the
coordinates $(x^{a})$ such that $\sum_{a=1}^{h}x^a=0$ (notation of (2.11)).
Accordingly, the orbits of $G_{\xi}$
are the $(h-1)$-dimensional planes through the points of $J^{-1}(\xi)$,
where the coordinates in the planes are
$(x^{a})$ with $\sum_{a=1}^{h}x^{a}=0$. Therefore, we have a quotient manifold
$$J^{-1}(\xi)/{\cal G}_{\xi}={\bf R}^{2(n-h)}=\{(x^{u},y_{u})_{u=h+1}^{n})\}.
\leqno{(2.15)}$$
The restriction of the form $\omega$ to $J^{-1}(\xi)$ is
$\omega'=\sum_{u=h+1}^{n}dx^{u}\wedge dy_{u}$, and this also is the expression
of the reduced symplectic form induced by $\omega$.
We see that we are in a situation where Proposition 2.1 and Theorem 1.1
apply, and the reduced
symplectic connection is again the natural flat torsionless connection of
${\bf R}^{2(n-h)}$.
\section{Presymplectic Manifolds}
From the geometrical point of view, the basic configuration where
symplectic reduction is encountered is that of a presymplectic manifold
i.e., a $(2n+p)$-dimensional differentiable manifold
$M$ endowed with a closed $2$-form $\omega$ of constant rank $2n$.
Then, the $p$-dimensional distribution
$$I:=span\{X\in \Gamma TM\:/\:i(X)\omega=0\}$$
is integrable, and yields a foliation ${\cal I}$, $T{\cal I}=I$, called the
{\em characteristic foliation}. Furthermore, if the space of leaves
$Q=M/{\cal I}$ is a Hausdorff manifold,
$Q$ is endowed with the {\em reduced symplectic form} $\omega'$
which is the projection of $\omega$, and $(Q,\omega')$ is the {\em reduced
symplectic manifold} of $(M,\omega)$ \cite{LM}.

Let us define a {\em presymplectic connection} on a
presymplectic manifold $(M,\omega)$ as being a linear connection $\nabla$
on $M$ which is $\omega$-compatible i.e., $\nabla\omega=0$, and its torsion
$T_{\nabla}$ takes values in the vector bundle $I$. Then, let us assume that
$(M,\omega)$ is endowed with a presymplectic connection $\nabla$.
We want to find conditions which ensure the existence of a
corresponding induced symplectic connection on the reduced manifold
$(Q,\omega')$. We begin with
\proclaim 3.1 Proposition. Let $V$ be a differentiable manifold endowed
with a foliation ${\cal F}$ and with a linear connection
$\nabla$ such that: i) ${\cal F}$ is parallel with respect to $\nabla$
(i.e., $\nabla(T{\cal F})\subseteq T{\cal F}$), ii) the torsion
$T_{\nabla}$ of $\nabla$ takes values in $T{\cal F}$,
iii) the curvature $R_{\nabla}$
of the connection $\nabla$ satisfies the condition
$$R_{\nabla}(Z,X)Y\in\Gamma T{\cal F},\hspace{5mm}
\forall Z\in\Gamma
T{\cal F},\,\forall X,Y\in\Gamma TV.\leqno{(3.1)}$$
Then $\nabla$ induces an ${\cal F}$-projectable connection on the normal
bundle $\nu{\cal F}$ of the foliation ${\cal F}$. \par
\noindent{\bf Proof.} By definition, $\nu{\cal F}=TV/T{\cal F}$, and a
cross section $\sigma\in\Gamma \nu{\cal F}$ is an equivalence class
$\sigma=[Y]_{\cal F}$ where $Y\in\Gamma TV$, and $Y_{1},Y_{2}$ yield the
same class $\sigma$ iff $Y_{1}-Y_{2}\in\Gamma T{\cal F}$.
Accordingly, the result of the formula
$$\nabla'_{X}\sigma:=[\nabla_{X}Y]_{\cal F}\hspace{1cm}(X\in\Gamma TV)
\leqno{(3.2)}$$
is independent of the choice of $Y$ (hypothesis i)), and $\nabla'$ is a
well defined connection on $\nu{\cal F}$.

We want to show that $\nabla'$ is ${\cal F}$-projectable i.e., it projects
to connections of the local {\em slice spaces} of ${\cal F}$. The
conditions for this are: a) if $\sigma$ is a projectable cross
section and $X\in\Gamma T{\cal F}$ then $\nabla'_{X}\sigma=0$,
b) if $\sigma$ is a projectable cross section
of $\nu{\cal F}$ and $X$ is a projectable vector
field then $\nabla'_{X}\sigma$ is a projectable cross section of
$\nu{\cal F}$ \cite{Ml}.
We recall that the projectability of $\sigma$ means that the
vector fields which represent $\sigma$ are projectable, which, in turn,
means that their flow preserves the foliation ${\cal F}$.

Condition a) easily follows from i) and ii).

For b), notice that $\forall Z\in\Gamma T{\cal F},\forall
X,Y\in\Gamma_{pr}TM$ one has
$$\nabla_{Z}\nabla_{X}Y=\nabla_{X}\nabla_{Z}Y+
\nabla_{[Z,X]}Y+R_{\nabla}(Z,X)Y\in \Gamma T{\cal F},$$
because of a) and of hypothesis iii). (In (3.1), we didn't have to ask
$X,Y\in\Gamma_{pr}TM$ since, $R_{\nabla}$ being a tensor, (3.1) is a
pointwise condition.) Accordingly, and using i) and ii) again, we get
$$[Z,\nabla_{X}Y]=\nabla_{Z}\nabla_{X}Y-\nabla_{\nabla_{X}Y}Z
-T_{\nabla}(Z,\nabla_{X}Y)\in
\Gamma T{\cal F},$$ therefore, $\nabla_{X}Y\in\Gamma_{pr}TM$. Q.e.d.

Using Proposition 3.1 we get
\proclaim 3.2 Proposition. Let $(M,\omega)$ be a presymplectic manifold
with the characteristic foliation ${\cal I}$, and let $\nabla$ be a
presymplectic connection on $M$. Assume that the curvature of $\nabla$
satisfies the condition
$$i(R_{\nabla}(Z,X)Y))\omega=0,\hspace{5mm}\forall Z\in\Gamma
T{\cal I},\,\forall X,Y\in\Gamma TM. \leqno{(3.3)}$$
Then, if the reduced symplectic manifold $(Q=M/{\cal I},\omega')$ exists, it
has a symplectic connection $\nabla'$ induced by $\nabla$. \par
\noindent{\bf Proof.} It is easy to see that $\nabla\omega=0$ implies
hypothesis i) of Proposition 3.1
for the foliation ${\cal I}$ on $M$. Hypothesis ii) holds because of the
definition of a presymplectic connection, and iii) follows from (3.3).
Therefore, there exists an
induced projectable connection $\nabla'$ on $\nu{\cal I}$,
which, obviously, projects to a connection on $Q$ that we also denote by
$\nabla'$. From the
definition of $\nabla'$ as given in the proof of Proposition 3.1,
and since $T_{\nabla}$ takes values in $I$, it follows
that, on $Q$, $\nabla'$ has no torsion and that
$\nabla\omega=0$ implies $\nabla'\omega'=0$. Q.e.d.

We remember that hypothesis iii) of Proposition 3.1
is what ensures the
projectability of the ``${\cal I}$-normal part" of $\nabla$. We want to
notice that this property may also be expressed without the use of the
normal bundle of ${\cal I}$. Namely, if $V$ is a manifold with a foliation
${\cal F}$, we will say that a linear connection $\nabla$ on $V$ is
${\cal F}$-{\em adapted} if $\forall X\in\Gamma T{\cal F}$ and $\forall
Y\in\Gamma_{pr}TV$, $\nabla_{X}Y\in \Gamma T{\cal F}$. Then if, moreover,
$\forall X,Y\in\Gamma_{pr}TV$ one has $\nabla_{X}Y\in\Gamma_{pr}TV$, we will
say that the connection $\nabla$ is ${\cal F}$-{\em projectable}.
These definitions are correct since they are
compatible with multiplication of a vector
field in $\Gamma{\cal F}$ by an arbitrary function, and with multiplication
of $X,Y\in\Gamma_{pr}TV$ by a projectable function. It is easy to see that
any foliated manifold has adapted connections, and, if the foliation has a
transversal projectable connection, the manifold has a projectable
connection.

On a presymplectic manifold $(M,\omega)$, any presymplectic connection
$\nabla$ is adapted to
the characteristic foliation ${\cal I}$. Indeed, ${\cal I}$
is $\nabla$-parallel, and the condition which was imposed on $T_{\nabla}$
yields $$\nabla_{X}Y=\nabla_{Y}X+[X,Y]+T_{\nabla}(X,Y)\in\Gamma T{\cal I},
\leqno{(3.4)}$$
$\forall X\in\Gamma T{\cal I},\,\forall Y\in\Gamma_{pr}TM$.
Using this remark, it follows
\proclaim 3.3 Proposition. If the reduced symplectic manifold $(Q=M/{\cal
I},\omega')$ of the presymplectic manifold $(M,\omega)$ exists, any
${\cal I}$-projectable presymplectic connection $\nabla$ on $M$ induces
a symplectic connection $\nabla'$ on $Q$. \par
\noindent{\bf Proof.} Any vector fields ${\cal X},{\cal Y}\in \Gamma TQ$ are
projections of vector fields $X,Y\in\Gamma_{pr}TM$, and if we put
$\nabla'_{\cal X}{\cal Y}=pr(\nabla_{X}Y)$ we are done. Q.e.d.

Proposition 3.3 raises the question of existence of
projectable presymplectic connections on a given presymplectic manifold
$(M,\omega)$.

First, we use the known technique of the almost symplectic
case e.g., \cite{V2} to get all the connections which satisfy
$\nabla\omega=0$. Let us fix a transversal distribution $S$ of the
characteristic distribution $I$ ($TM=S\oplus I$), a connection $D^{I}$ on
the vector bundle $I$, and a connection $\stackrel{\circ}{D^S}$ on the
vector bundle $S$. Then, change the connection
$\stackrel{\circ}{D^S}$ to the Bott connection \cite{Ml}
$$D^{S}_{(X'+X'')}Y'=pr_{S}[X'',Y']+\stackrel{\circ}{D^S}_{X'}Y',
\leqno{(3.5)}$$
where $X',Y'\in\Gamma S,\,X''\in\Gamma I$, and define the connection
$$D=D^{S}\oplus D^{I}\leqno{(3.6)}$$ on $TM$. Furthermore, let us define
$\Theta\in\Gamma End(TM\wedge TM,TM)$ by
$$\begin{array}{l}
\Theta(X'',Y'')=0,\;\Theta(X',Y'')=0,\;\Theta(X',Y')\in\Gamma
S,\vspace{2mm}\\ \omega(\Theta(X',Y'),Z')=\frac{1}{2}
(D_{X'}\omega)(Y',Z'), \end{array} \leqno{(3.7)}$$
$\forall X',Y'\in \Gamma S,\,\forall X'',Y''\in\Gamma I$.
(The result is well defined because $\omega$ is non degenerate on $S$.)

Using $\Theta$, we get a new connection on $M$ namely,
$$\stackrel{\circ}{\nabla}=D+\Theta,\leqno{(3.8)}$$
and the evaluation of $(\stackrel{\circ}{\nabla}\omega)(X,Y)$ for
$X,Y\in\Gamma_{pr}TM$ yields $\stackrel{\circ}{\nabla}\omega=0$. Hence, we
obtained one $\omega$-compatible connection. Necessarily, all the others
are given by
$$\nabla=\stackrel{\circ}{\nabla}+A,\leqno{(3.9)}$$
where $A\in\Gamma End(TM\wedge TM,TM)$ satisfies
$$\omega(A(X,Y),Z)+\omega(Y,A(X,Z))=0.\leqno{(3.10)}$$

The torsion of the connection $D$ given by (3.6) is
$$T_{D}(X'+X'',Y'+Y'')= \stackrel{\circ}{D^S}_{X'}Y'
-\stackrel{\circ}{D^S}_{Y'}X'\leqno{(3.11)}$$
$$-pr_{S}[X',Y']+\:{\rm
term\; in }\:\Gamma I$$ (notation of (3.5)).
Hence, if we take $\stackrel{\circ}{D}=pr_{S}\circ K$, where $K$
is a torsionless covariant derivative on $M$, the associated connection $D$
will have an $I$-valued torsion, and it will follow
$$d\omega(X,Y,Z)=\sum_{cycl.(X,Y,Z)}(D_{X}\omega)(Y,Z)=0.\leqno{(3.12)}$$

Furthermore, the torsion of the connection (3.9) is
$$T_{\nabla}(X,Y)=T_{D}(X,Y)+\Theta (X,Y)-\Theta (Y,X)+A(X,Y)-A(Y,X),
\leqno{(3.13)}$$
where $X,Y\in\Gamma TM$, $T_{D}(X,Y)\in\Gamma I$,
and from (3.7), (3.12) we obtain
$$\omega(T_{\nabla}(X',Y'),Z')=\omega(A(X',Y'),Z')-\omega(A(Y',X'),Z')
\leqno{(3.14)}$$
$$-\frac{1}{2}(D_{Z'}\omega)(X',Y'),$$
where $X',Y',Z'\in\Gamma TS$.
Then, if we ask $A(X'',Y)=0$ for $X''\in\Gamma I,\,
Y\in\Gamma TM$, and $A(X',Y')\in\Gamma S$, and such that \cite{V2}
$$\omega(A(X',Y'),Z')=\frac{1}{6}\{(D_{Y'}\omega)(X',Z')
+(D_{Z'}\omega)(X',Y')\},\leqno{(3.15)}$$
we obtain a connection $\nabla$ which satisfies the condition
$$\omega(T_{\nabla}(X',Y'),Z')=0,$$
therefore, $\nabla$ has an $I$-valued torsion.

With a closer look at this latter connection we see that, in fact, the
following result holds
\proclaim 3.4 Proposition. Every presymplectic manifold $(M,\omega)$ has
presymplectic connections. Moreover, the manifold has a
projectable presymplectic connection
with respect to the characteristic foliation ${\cal I}$, iff the normal
bundle of ${\cal I}$ has a projectable connection. \par
\noindent{\bf Proof.} The first part was already proven above. For the
second part, if the normal bundle $\nu{\cal I}$ has a projectable
connection, we may use it as $\stackrel{\circ}{D^S}$ of (3.5), and we see
that the connection $\nabla$ defined by the use of (3.15) is projectable.
Conversely, if $\nabla$ is a projectable presymplectic connection on
$(M,\omega)$, and if we put $$D_{X}[Y]_{\cal I}=[\nabla_{X}Y]_{\cal I}
\hspace{5mm}\forall X,Y\in\Gamma TM$$
we get a projectable connection on $\nu{\cal I}$. Q.e.d.

We recall that a foliation has a projectable connection on its normal
bundle iff
the {\em Atiyah class} of the foliation vanishes \cite{Ml}.
 \vspace{1cm}
{\small Department of Mathematics, \\}
{\small University of Haifa, Israel.\\}
{\small E-mail: vaisman@math.haifa.ac.il}
\end{document}